\documentclass[final]{elsarticle}

\usepackage{amsfonts}
\usepackage{enumerate} 
\usepackage{amsthm}
\usepackage{amsmath}
\usepackage{mathtools}
\usepackage{amscd}
\usepackage[active]{srcltx} 
\usepackage{latexsym}
\usepackage{amssymb}
\usepackage{enumitem}

\newtheorem{thm}{Theorem}[section]
\newtheorem{prop}[thm]{Proposition}
\newtheorem{cor}[thm]{Corollary}
\newtheorem*{cor*}{Corollary}
\newtheorem{lema}[thm]{Lemma}
\newtheorem*{lema*}{Lemma}

\numberwithin{equation}{section}
\theoremstyle{definition}
\newtheorem*{Def}{Definition}

\newenvironment{dem}{\vspace{1ex}\noindent{\it Proof.}\hspace{0.5em}}
{\hfill\qed\vspace{1ex}}

\newtheorem{example}{Example}

\newtheorem*{obs*}{Remark}
\newtheorem*{thm*}{Theorem}
\newtheorem*{prop*}{Proposition}

\newtheoremstyle{dotless}{}{}{}{}{}{}{ }{}
\theoremstyle{dotless}

\newcommand{\matriz}[4]{\displaystyle\
    \left(
       \begin{array}{cc}
        {#1}&{#2}\\
        {#3}&{#4}
       \end{array}
     \right)}

\newcommand{\PI}[2]{\left\langle \,#1 , #2\, \right\rangle}
\newcommand{\PIW}[2]{\left\langle \,#1 , #2\, \right\rangle_W}

\newcommand{\NW}[1]{\Vert #1 \Vert_{W}^{2}}
\newcommand{\NC}[1]{\Vert #1 \Vert}

\newcommand{\NWW}[1]{\Vert #1 \Vert_{W}}

\newcommand{\Nphis}[1]{\Vert #1 \Vert_{p}}

\newcommand{\Nphiss}[1]{\Vert #1 \Vert_{p}}
\newcommand{\Nphissw}[1]{\Vert #1 \Vert_{p,W}}

\newcommand{\ra}{\rightarrow}

\newcommand{\St}{\mathcal{S}}
\newcommand{\Rt}{\mathcal{R}}
\newcommand{\HH}{\mathcal{H}}

\newcommand{\M}{\mathcal{M}}

\newcommand{\CC}{\mathbb{C}}
\newcommand{\Q}{\mathcal{Q}}

\newcommand{\noi}{\noindent}
\newcommand{\dit}[1]{\textit{#1}}

\journal{Linear Algebra and its Applications}

\begin{document}

\begin{frontmatter}

\title{Weighted least squares solutions of the equation $AXB-C=0$}
\author[UBA]{Maximiliano Contino}
\ead{mcontino@fi.uba.ar}

\author[IAM,UBA]{Juan Giribet}
\ead{jgiribet@fi.uba.ar}

\author[IAM,UBA]{Alejandra Maestripieri \corref{ca}}
\ead{amaestri@fi.uba.ar}

\address[IAM]{Instituto Argentino de Matem\'atica ``Alberto P. Calder\'on''\\ Saavedra 15, Piso 3 (1083) Buenos Aires, Argentina}
\address[UBA]{Departamento de Matem\'atica-- Facultad de Ingenier\'{\i}a -- Universidad de Buenos Aires\\ Paseo Col\'on 850 (1063) Buenos Aires, Argentina}

\cortext[ca]{Corresponding author}

\begin{abstract}
		Let $\HH$ be a Hilbert space, $L(\HH)$ the algebra of bounded linear operators on $\HH$ and $W \in L(\HH)$ a positive operator such that $W^{1/2}$ is in the p-Schatten class, for some $1 \leq p< \infty.$ Given $A, B \in L(\HH)$ with closed range and $C \in L(\HH),$ we study the following weighted approximation problem: analize the existence of 
\begin{equation}
		\underset{X \in L(\HH)}{min}\Nphissw{AXB-C}, \label{eqa1}
		\end{equation}
		where $\Nphissw{X}=\Nphiss{W^{1/2}X}.$ 
		We also study the related operator approximation problem: analize the existence of 
		\begin{equation}
		\underset{X \in L(\HH)}{min} (AXB-C)^{*}W(AXB-C), \label{eqa2}
		\end{equation}
	    where the order is the one induced in $L(\HH)$ by the cone of positive operators.
		In this paper we prove that the existence of the minimum of \eqref{eqa2} is equivalent to the existence of a solution of the normal equation $A^*W(AXB-C)=0.$ We also give sufficient conditions for the existence of the minimum of \eqref{eqa1} and we characterize the operators where the minimum is attained.
		
\end{abstract}

\end{frontmatter}

\section{Introduction}
In signal processing language, {\it sampling} is an operation which converts a continuous signal (modelled as a vector in an adecuate Hilbert space $\HH$) into a discrete one.
  
Frequently the samples of a signal $f\in\HH$ are represented in the the following way: given a frame $\{v_n\}_{n\in\mathbb{N}}\subseteq\HH$ of a closed subspace $\St$, called the sampling subspace, the samples are given by $\{f_n\}_{n\in\mathbb{N}}=\{\PI{f}{v_n}\}_{n\in\mathbb{N}}\in\ell^2(\mathbb{N})$.
On the other hand, given samples $\{f_n\}_{n\in\mathbb{N}}\in\ell^2(\mathbb{N})$, the reconstructed signal $\hat{f}$ is given by $\hat{f}=\sum_{n\in\mathbb{N}} f_n w_n$, where $\{w_n\}_{n\in\mathbb{N}}$ is a frame of the closed subspace $\Rt$,  called the reconstruction subspace. 

Suppose that, $A$ and $B$ are the {\it synthesis operators} corresponding to the frames $\{w_n\}_{n\in\mathbb{N}}$ and $\{v_n\}_{n\in\mathbb{N}}$, respectively, i.e. $A, B: \ell^2 (\mathbb{N})\rightarrow\HH$ are the operators such that, if $x=\{x_n\}_{n\in\mathbb{N}}\in\ell^2 (\mathbb{N})$, $Ax=\sum_{n\in\mathbb{N}} x_n w_n$ and $Bx=\sum_{n\in\mathbb{N}} x_n v_n$, which are bounded since $\{v_n\}_{n\in\mathbb{N}}$ and $\{w_n\}_{n\in\mathbb{N}}$ are frames.
Observe that, the samples of $f$ are given by $\{f_n\}_{n\in\mathbb{N}}=B^*f$ and, given samples $\{f_n\}_{n\in\mathbb{N}}$ the reconstructed signal is given by $\hat{f}=A(\{f_n\}_{n\in\mathbb{N}})$, see \cite{Eldar}, \cite{Unser}.

If we only know the samples of a signal $\{f_n\}_{n\in\mathbb{N}}\in\ell^2(\mathbb{N})$, in general it is not possible to recover the signal $f\in\HH$, even if we apply a digital filter (a bounded linear operator $X:\ell^2(\mathbb{N})\rightarrow\ell^2(\mathbb{N})$) to these samples. But, in some cases it is possible to find a good representation of the signal $f\in\HH$, i.e., a recovered signal $\hat{f}=AXB^*f$ that has good properties. For instance, in the classical sampling scheme (where sampling and reconstruction subspaces coincide) it is possible to reconstruct the best approximation of the signal $f$, i.e., it is possible to find $X$ such that $AXB^*=P_\St$ and then $\hat{f}=P_\St f$, where $P_\St$ is the orthogonal proyection onto $\St=\Rt=R(A)$. Another interesting example, where the sampling and reconstruction subspaces may not coincide,  is the so called consistent sampling scheme, where the samples of the reconstructed signal $\hat{f}$ are equal to the samples of the original signal $f$, i.e. $B^*\hat{f}=B^*f$, in this case $X$ is such that $Q=AXB^*$ turns out to be an oblique projection. Consequently, the reconstructed signal $\hat{f}$ is not necessarily a good approximation of $f$, since the distance $\|f-\hat{f}\|=\|f-AXB^*f\|$ is not minimized. Now suppose we want to find a digital filter $X\in\ell^2(\mathbb{N})\rightarrow\ell^2(\mathbb{N})$ such that $AXB^*f$ is a good approximation of $f$ in $R(A)=\Rt$, i.e. we want that $AXB^*$ approximates $P_\Rt$ in some sense. For instance, we may want to find $X_0:\ell^2(\mathbb{N}) \rightarrow \ell^2(\mathbb{N})$ a bounded linear operator, such that, for every $f\in\HH$
$$
\|(AX_0B^*-P_\Rt)f\|\leq\|(AXB^*-P_\Rt)f\|,
$$
for every $X\in L(\HH)$ (the algebra of linear bounded operators on $\HH$). This means that, we are interested in the following problem,
$$
\min_{X\in L(\ell^2(\mathbb{N}))} (AXB^*-P_\Rt)^*(AXB^*-P_\Rt),
$$
with the order induced in $L(\HH)$ by the cone of positive operators.

Alternatively, we can approximate $P_\Rt$ in some convenient operator norm. For example, in the finite dimensional setting, it is usual to consider the Frobenius norm; the associated problem becomes studying the existence of 
$$
\min_{X\in  L(\ell^2(\mathbb{N}))} \|AXB^*-P_\Rt\|_2,
$$

where $\| \cdot \|_2$ is the Frobenius norm.

In this work we are interested in studying an extension of these problems. More specifically, 
given $A, B \in L(\HH)$ with closed range, $C \in L(\HH),$ $W \in L(\HH)$ positive, we study the existence of 
\begin{equation}
\underset{X \in L(\HH)}{min} (AXB-C)^{*}W(AXB-C), \label{eq122}
\end{equation}
with the order induced in $L(\HH)$ by the cone of positive operators.
If $W$ also satisfies that $W^{1/2} \in S_p,$ the p-Schatten class (for some $p$ with $1 \leq p < \infty$), consider the seminorm associated to $W$, $$\Nphissw{X}=\Nphiss{W^{1/2}X},$$
for $X\in L(\HH).$
We study the existence of
\begin{equation} 
\underset{X \in L(\HH)}{min}\Nphissw{AXB-C}. \label{eq131}
\end{equation}
We are also interested in giving a characterization of the set of solutions of these problems.

\vspace{0,3cm}
There are several examples of these minimization problems, in \cite{CMR2014} a similar problem related with frame theory is studied,
in \cite{Gold1} the existence of minimum of $\NC{AX-I}_{p}$ in the finite dimensional setting is given. In \cite{Mah2}, \cite{Mah3}, \cite{Mecheri},  \cite{Mecheri2} and in \cite{Bounkhel}, the existence of minimum of $\NC{AX-C}_{p}$ in Hilbert spaces,  with suitable hypotesis to warrant that $AX-C \in S_p$, was studied using differentiation techniques and also in \cite{Nashed}, where a connection between $p$-Schatten norms and the order in $L(\HH)^+$ (the cone of semidefinite positive operators) is established. In \cite{Contino}, the existence of minimum of $\NC{AX-C}_{p,W}$ in Hilbert spaces was stated, however  the introduction of an operator $B \in L(\HH)$ with closed range, produces notable differences in the final results. 

The existence of minimum of $\NC{AXB-C}_{p}$ in Hilbert spaces, with suitable hypotesis to warrant that $AXB-C \in S_p$, was studied in \cite{Mah4} and in \cite{Mah3} using differentiation techniques. In \cite{Changsen}, a characterization of the critical points (and equivalently the global minima) of the map $\NC{AXB-C}^p_{p}$ is given.
However, the introduction of a weight $W\in L(\HH)^+$ plays an important role, since we are introducing on $\HH$ a semi-inner product associated to $W$ for which $\HH$ is no longer a Hilbert space, unless $W$ is invertible.

\vspace{0,3cm}
The contents of the paper are the following. In section 2, the concept of $W$-inverse of an operator $A$ in the range of an operator $B$ is introduced, together with some properties. Some results of shorted operators and  compressions are stated. Also, definitions and properties of directional derivatives are included.

In section 3, we study problem \eqref{eq122}. We prove that if $N(B) \subseteq N(A^*WC)$ then the infimum of the set $\{(AXB-C)^{*}W(AXB-C): X \in L(\HH)\}$ (where the order is the one induced by the cone of positive operators) exists and it is equal to $C^*W_{/R(A)}C,$ where $W_{/R(A)}$ is the shorted operator of $W$ to $R(A).$ We also prove that the minimum of the previous set exists if and only if $N(B)\subseteq N(A^*WC)$ and $R(C) \subseteq R(A)+W(R(A))^{\perp}.$ Moreover, we prove that an operator $X_0$ minimize this problem, if and only if $X_0B$ is a $W$-inverse of $A$ in $R(C).$ 

 In section 4, it is shown that if $W^{1/2}$ is in the $p$-Schatten class, for some $1 \leq p< \infty,$ then if $N(B)\subseteq N(A^*WC)$ and $R(C) \subseteq R(A)+W(R(A))^{\perp},$ then the minimum of the set $\{ \Nphiss{W^{1/2}(AXB-C)} : X \in L(\HH)\}$ exists. In Lemma \ref{critical_point}, we give a characterization of the critical points (and equivalently the global minima) of the map $\NC{AXB-C}^p_{p,W}$, which is similar to the one considered in \cite{Changsen}, with the introduction of a weight $W$ such that $W^{1/2}$ is in the $p$-Schatten class, for some $1 \leq p< \infty.$ If $p=2$ or alternatively $1 \leq p < \infty$ and $N(B) \subseteq N(A^*WC),$ it is proven that the existence of the minimum of the previous set is equivalent to the existence of solution of the normal equation $A^*W(AXB-C)B^*=0.$ Finally, some examples are given to show that, in general, the existence of the minimum of the previous set is not equivalent to the existence of the solution of the presented normal equation, showing that in \cite[Theorem~4.1]{Mah4} additional hypothesis should be added. 

\section{Preliminaries}
Throughout $\HH$ denotes a separable complex Hilbert space, $L(\HH)$ is the algebra of bounded linear operators from $\HH$ to $\HH$, and $L(\HH)^{+}$ the cone of semidefinite positive operators. $GL(\HH)$ is the group of invertible operators in $L(\HH),$ $CR(\HH)$ is the subset of $L(\HH)$ of all operators with closed range.
For any $A \in L(\HH),$ the range and the nullspace of $A$ are denoted by $R(A)$ and $N(A)$ respectively. Finally,  $A^{\dagger}$ denotes the  Moore-Penrose inverse of the operator $A \in L(\HH)$. For $A, B \in L(\HH)^{+},$ $A \leq B$ if $B-A \in L(\HH)^+.$

Given a closed subspace $\M,$  $P_{\M}$ denotes the orthogonal projection onto $\M.$ Also, $\Q$ denotes the subset of $L(\HH)$ of oblique projections, i.e. $\Q=\{Q \in L(\HH): Q^{2}=Q\}.$ 

Given $W\in L(\HH)^{+},$ $\PIW{x}{y}=\PI{Wx}{y}, \ x, y \in \HH$ defines a semi-inner product on $\HH.$ There is also a seminorm associated to $W,$ namely $\NW{x}=\PI{Wx}{x}, \ x \in \HH.$

The $W$-orthogonal complement of  $\St\subseteq \HH$ is given by $$\St^{\perp_{W}}=\{x \in \HH: \PI{Wx}{y}=0, \ y\in \St\}=W^{-1}(\St^{\perp}).$$

\vspace{0,3cm}
We now give the definitions of $W$-least squares solution of the equation $Az=x.$

\begin{Def} Given $A\in CR(\HH),$ $W\in L(\HH)^{+}$ and $x\in \HH,$ $u\in \HH$ is a $W$-least squares solution or $W$-$LSS$ of $Az=x,$ if $$\NWW{Au-x}\leq\NWW{Az-x}, \mbox{ for every } z \in \HH.$$
\end{Def}

The next theorem describes some properties of the $W$-least squares solutions of $Az=x.$

\begin{thm}
\label{thmWLSS} Given $A\in CR(\HH),$ $W\in L(\HH)^{+}$ and  $x \in \HH,$ there exists a $W$-$LSS$ of $Az =x$ if and only if $x \in R(A)+R(A)^{\perp_W}.$ 
\end{thm}

\begin{dem} See \cite{Spline}.\\
\end{dem}

\vspace{0,3cm}
In \cite{Mitra} S. K. Mitra and C. R. Rao introduced the notion of the $W$-inverse of a matrix. Observe that, in this setting it holds $\HH=R(A)+R(A)^{\perp_W}$, because $\HH$ is a finite dimensional space, \cite{Shorted3}. This concept was extended to operators in \cite{WGI} and in \cite{Contino}.  

\begin{Def} Given $A \in CR(\HH),$ $B \in L(\HH)$ and $W\in L(\HH)^{+},$ $X_0 \in L(\HH)$ is a $W$-inverse of $A$ in $R(B),$ if for each $x \in \HH$, $X_0x$ is a $W$-$LSS$ of $Az=Bx,$ i.e. $$\NWW{AX_0x-Bx}\leq\NWW{Az-Bx}, \mbox{ for every } x, z \in \HH.$$  
When $B=I,$ $X_0$ is called the $W$-inverse of $A.$ See \cite{WGI}.
\end{Def}

\begin{thm}
\label{thmWinversa} Given $A\in CR(\HH), B \in L(\HH)$ and $W\in L(\HH)^{+},$ the following conditions are equivalent:
\begin{enumerate} 
\item [i)] The operator $A$ admits a $W$-inverse in $R(B),$
\item [ii)] $R(B) \subseteq R(A) + R(A)^{\perp_{W}},$ 
\item [iii)] the normal equation $A^{*}W(AX-B)=0$ admits a solution. 
\end{enumerate}
\end{thm}

\begin{dem}
\cite[Theo~2.4]{Contino}.\\
\end{dem}

\begin{cor} \label{Winversa 3} If $R(B) \subseteq R(A) + R(A)^{\perp_{W}},$ then the set of $W$-inverses of $A$ in $R(B)$ is the set of solutions of the equation $A^{*}W(AX-B)=0$, or equivalently the affine manifold
$$(A^{*}WA)^{\dagger}A^{*}WB+ \{L \in L(\HH): R(L) \subseteq N(A^{*}WA)\}.$$
\end{cor}

\vspace{0,3cm}
Given $W\in L(\HH)^{+}$ and a closed subspace $\St \subseteq \HH$ the notion of shorted operator of $W$ to $\St,$ was introduced by M. G. Krein in \cite{Krein} and later rediscovered by W. N. Anderson and G. E. Trapp who proved in \cite{Shorted2}, that the set
$\{ X \in L(\HH): \ 0\leq X\leq W \mbox{ and } R(X)\subseteq \St^{\perp}\}$ has a maximum element.

\begin{Def} The shorted operator of $W$ to $\St$ is defined by
$$W_{/\St}=\mbox{max } \{ X \in L(\HH): \ 0\leq X\leq W \mbox{ and } R(X)\subseteq \St^{\perp}\}.$$

The $\St$-compression $W_{\St}$ of $W$ is defined by $$W_{\St}=W-W_{/\St}.$$
\end{Def}

For many results on the notions of shorted operators, the reader is referred to \cite{Shorted1} and \cite{Shorted2}. 

Next we collect some results regarding $W_{/\St}$ and $W_{\St}$ which are relevant in this paper.

\begin{thm} \label{TeoShorted}
Let $W\in L(\HH)^{+}$ and $\St \subseteq \HH$ a closed subspace. Then  
\begin{enumerate} 
\item [i)] $W_{/\St}=\mbox{ inf } \{ E^{*}WE: E^2=E, \ N(E)=\St\};$ in general, the infimum is not attained,
\item [ii)] $R(W) \cap \St^{\perp} \subseteq R(W_{/\St}) \subseteq R(W^{1/2}) \cap \St^{\perp},$
\item [iii)] $\overline{N(W)+\St} \subseteq N(W_{/\St})=W^{-1/2}(\overline{W^{1/2}(\St)}),$
\item [iv)] $N(W_{\St})=W^{-1}(\St^{\perp})$ and $W(\St) \subseteq R(W_{\St}) \subseteq \overline{W(\St)}.$
\end{enumerate}
\end{thm}
The reader is referred to \cite{Shorted2} and \cite{Shorted3} for the proof of these facts.

\vspace{0,3cm}
\begin{Def} Let $T\in L(\HH)$ be a compact operator. By $\{\lambda_k(T)\}_{k\geq1}$ we denote the eigenvalues of $\vert T \vert = (T^{*}T)^{1/2},$ where each eigenvalue is repeated according to its multiplicity. Let $1\leq p < \infty,$ we say that $T$ belongs to the p-Schatten class $S_p,$ 
if $$\sum_{k\geq1}^{} \lambda_{k}(T) ^{p}<\infty,$$ and we note $$\Nphiss{T}= (\sum_{k\geq1}^{} \lambda_{k}(T) ^{p})^{1/p},$$ where $\Nphiss{\cdot}$ is called the p-Schatten norm.
\end{Def}
 
The reader is referred to \cite{Ringrose, Simon} for further details.\\

\begin{prop} \label{Prop Nashed} Let $1\leq p < \infty$, $T \in L(\HH)$ and $S\in S_p$. If  $T^{*}T\leq S^{*}S$ then $\| T \|_p \leq \| S \|_p$.
\end{prop}
 
\begin{dem} \cite[Prop~2.9]{Contino}. See also \cite[Prop~2.5]{Nashed}, where a more general result is given.
\end{dem}

\vspace{0,3cm}
The following theorem characterized the existence of solution of the equation $AXB=C.$

\begin{thm} \label{Teosol} Let $A, B, C \in L(\HH).$ If $R(A), R(B)$ or $R(C)$ is closed, then the equation $AXB=C$ admits a solution if and only if $R(C) \subseteq R(A)$ and $R(C^*) \subseteq R(B^*).$ 

In this case, the general solution of the equation $AXB-C=0$ is
$$A^{\dagger}CB^{\dagger}+L-A^{\dagger}ALBB^{\dagger},$$ for arbitrary $L \in L(\HH).$
\end{thm}

\begin{dem}
	See \cite[Theo.~3.1]{Arias}.
\end{dem}

\vspace{0,5cm}
Finally, we give a definition for the derivative of a real-valued function on a Banach space, that will be instrumental to prove some results of this paper.

\begin{Def} Let $(\mathcal{E}, \NC{\cdot})$ be a Banach space and $f: \mathcal{E} \rightarrow \mathbb{R}.$ Let $\phi \in [0, 2\pi)$ and $h>0$, then the $\phi-$directional derivative of $f$ at a point $x\in \mathcal{E}$ in direction $y \in \mathcal{E}$ is defined by
$$ D_{\phi}f(x,y)=lim_{h\rightarrow 0^{+}} \frac{f(x+h e^{i\phi}y)-f(x)}{h}.$$
\end{Def}

Observe that if $f: \mathcal{E} \ra \mathbb{R},$ $f(x)=\Vert x \Vert,$ then $D_{\phi}f(x,y)$ is a subadditive, positive functional on $\mathcal{E},$ such that 
$$\vert D_{\phi}f(x,y) \vert \leq \Vert y \Vert, \mbox{ for every } x, y \in \mathcal{E}.$$ See  \cite[Prop.~1.3]{Drago}.

\begin{thm} \label{TeoD} Let $G_p: S_p \rightarrow \mathbb{R}^{+},$ $1\leq p<\infty,$ $G_p(X)=\Nphis{X}^{p},$ and let $X, Y \in S_p.$ Then, for all $\phi \in [0, 2\pi),$
\begin{itemize}
\item [i)] for $1< p<\infty,$ $G_p$ has a $\phi-directional$ derivative given by 
$$D_{\phi}G_p(X,Y)=p \ Re \ [e^{i\phi} tr ( \vert X \vert^{p-1} U^{*}Y)],$$ 
\item[ii)]  for $p=1,$ $G_1$ has a $\phi-directional$ derivative given by 
$$D_{\phi}G_1(X,Y)=Re \ [e^{i\phi} tr(U^{*}Y)] + \Vert P_{N(X^{*})}YP_{N(X)}\Vert_{1},$$

\end{itemize}
where $Re (z)$ is the real part of a complex number $z$, $tr(T)$ denotes the trace of the operator $T$ and $X=U\vert X\vert,$ is the polar decomposition of the operator $X,$ with $U$ the partial isometry such that $N(U)=N(X).$
\end{thm}

\begin{dem} See \cite[Theorem 2.1]{Aiken} and \cite[Theorem 2.1]{Drago}. 
\end{dem}

\begin{lema} \label{LemaM2} Let $(\mathcal{E}, \NC{\cdot})$ be a Banach space and $f: \mathcal{E} \rightarrow \mathbb{R},$ such that $f$ has a $\phi-directional$ derivative for every $\phi \in [0, 2\pi),$ at every point $x \in \mathcal{E}$ and in every direction $y \in \mathcal{E}$.  
If $f$ has a global minimum at $x_0 \in \mathcal{E},$ then 
$$\underset{0 \leq \phi < 2\pi}{inf} (D_{\phi}f(x_0,y)) \geq 0, \mbox{ for every } y \in \mathcal{E}.$$\end{lema}

\begin{dem} See \cite[Theorem 2.1]{Mecheri}.\end{dem}

\bigskip

\section{Minimization results in the operator order}
In this section we study the first problem mentioned in the introduction: given $A, B \in CR(\HH), \ C \in L(\HH), W \in L(\HH)^{+},$ we analize the existence of 
\begin{equation}
\underset{X \in L(\HH)}{inf} (AXB-C)^{*}W(AXB-C), \label{eq22}
\end{equation}
with the order induced in $L(\HH)$ by the cone of positive operators.

It was proven in \cite[Prop.~4.2]{Contino}, that if $A \in CR(\HH), \ C \in L(\HH)$ and $W \in L(\HH)^{+},$ then 
$$\underset{X \in L(\HH)}{inf} (AX-C)^*W(AX-C)=C^*W_{/R(A)}C.$$ Therefore,
\begin{equation}
(AXB-C)^{*}W(AXB-C) \geq C^*W_{/R(A)}C, \mbox{ for every } X \in L(\HH). \label{eqintro}
\end{equation}
The next result provides a sufficient condition for the existence of the infimum in \eqref{eq22}.
From now on, consider $$H(X)=(AXB-C)^{*}W(AXB-C).$$
\begin{prop} \label{Propinf2} Let $A, B \in CR(\HH), \ C \in L(\HH)$ and $W \in L(\HH)^{+}.$ 
If $N(B)\subseteq N(A^*WC)$ then the infimum of the set $\{H(X): X \in L(\HH)\}$  exists and
	$$\underset{X \in L(\HH)}{inf}H(X)=C^*W_{/R(A)}C.$$
\end{prop}

\begin{dem}
Suppose $N(B)\subseteq N(A^*WC).$ Then, it can be checked that
\begin{equation}
H(X)=G(X)+C^*WC-P_{N(B)^{\perp}}C^*WCP_{N(B)^{\perp}}, \label{eq23}
\end{equation}
where $G(X)=(AXB-CP_{N(B)^{\perp}})^{*}W(AXB-CP_{N(B)^{\perp}}).$

The set $\{G(X): X \in L(\HH)\}$ always admits an infimum. In fact,	let $X\in L(\HH), $ writing $W=W_{/R(A)}+W_{R(A)},$ it follows that
\begin{align*}
G(X)&=(AXB-CP_{N(B)^{\perp}})^{*}W(AXB-CP_{N(B)^{\perp}})\\
&=P_{N(B)^{\perp}}C^*W_{/R(A)}CP_{N(B)^{\perp}}+(AXB-CP_{N(B)^{\perp}})^{*}W_{R(A)}(AXB-CP_{N(B)^{\perp}})\\
&\geq P_{N(B)^{\perp}}C^*W_{/R(A)}CP_{N(B)^{\perp}},
\end{align*}
because $R(A) \subseteq N(W_{/R(A)})$ (see Theorem \ref{TeoShorted}). 
Hence $P_{N(B)^{\perp}}C^*W_{/R(A)}CP_{N(B)^{\perp}}$ is a lower bound of $G(X).$

If $D\geq 0$ is any other lower bound of $G(X),$ then $$D\leq G(X),\mbox { for every } X \in L(\HH).$$
In particular, $$D\leq P_{N(B)^{\perp}}C^*E^{*}WECP_{N(B)^{\perp}},$$ 
where $E$ is any projection such that $N(E)=R(A)$. In fact $R((I-E)CP_{N(B)^{\perp}})\subseteq R(I-E) =N(E)=R(A)$ and $R(P_{N(B)^{\perp}}C^*(I-E)^*)\subseteq R(P_{N(B)^{\perp}})=R(B^*)$, then by Theorem \ref{Teosol}, there exists $X_0 \in L(\HH),$ such that $(I-E)CP_{N(B)^{\perp}}=AX_0B,$ i.e., $(-E)CP_{N(B)^{\perp}}=AX_0B-CP_{N(B)^{\perp}}.$
Therefore, by \cite[Lemma.~4.1]{Contino} 
$$D \leq inf \{P_{N(B)^{\perp}}C^*E^{*}WECP_{N(B)^{\perp}} : \ E^{2}=E, \ N(E)=R(A) \}= P_{N(B)^{\perp}}C^*W_{/R(A)}CP_{N(B)^{\perp}}.$$ 
Thus, $$P_{N(B)^{\perp}}C^*W_{/R(A)}CP_{N(B)^{\perp}}=\underset{X \in L(\HH)}{inf}G(X).$$
Then, it follows that the infimum of $H(X)$ exists, moreover
\begin{align*}
\underset{X \in L(\HH)}{inf}H(X)&=\underset{X \in L(\HH)}{inf}G(X)+C^*WC-P_{N(B)^{\perp}}C^*WCP_{N(B)^{\perp}}\\
&=P_{N(B)^{\perp}}C^*W_{/R(A)}CP_{N(B)^{\perp}}+C^*WC-P_{N(B)^{\perp}}C^*WCP_{N(B)^{\perp}}\\
&=C^*WC-P_{N(B)^{\perp}}C^*W_{R(A)}CP_{N(B)^{\perp}}.
\end{align*}

But since $N(B) \subseteq N(A^*WC)$ and $N(W_{R(A)})=N(A^*W)$ (see Theorem \ref{TeoShorted}), we have
$W_{R(A)}CP_{N(B)}=0.$
	
Therefore
$$ \underset{X \in L(\HH)}{inf}H(X)=C^*WC-P_{N(B)^{\perp}}C^*W_{R(A)}CP_{N(B)^{\perp}} =C^*WC-C^*W_{R(A)}C=C^*W_{/R(A)}C.$$
\end{dem}

\vspace{0,3cm}

Now we give an example where Problem \ref{eq22} admits an infimum, but $N(B) \not \subseteq N(A^*WC),$ showing that the condition in Proposition \ref{Propinf2} is not necessary for the existence of infimum in \eqref{eq22}.

\example \label{example1} Let $\HH=\mathbb{C}^{2},$ $W=I,$ $A=C=\begin{bmatrix}
1 & 0 \\ 
0 & 0 \\
\end{bmatrix}$ and $B=\begin{bmatrix}
0 & 1 \\ 
0& 0 \\ \end{bmatrix}.$
Observe that  $N(B) \not \subseteq N(A^*WC).$

Let $X=\begin{bmatrix}
x & y \\ 
z & w \\
\end{bmatrix} \in \mathbb{C}^{2\times 2},$ where $x, y, z, w \in \mathbb{C}.$ Then

$AXB-C=\begin{bmatrix}
-1 & x \\ 
0 & 0 \\
\end{bmatrix},$ and
$H(X)=(AXB-C)^*(AXB-C)=\begin{bmatrix}
1 & -x \\ 
-\overline{x} & |x|^2 \\
\end{bmatrix}.$

Let $u,v \in \mathbb{C}$ then it can be checked that

$$\langle (AXB-C)^*(AXB-C) (u, v), (u,v) \rangle= |u-xv|^2.$$   
Since for any $u,v\in\CC$ there exists $x\in\CC$ such that $u-xv=0$, it follows that $$\underset{X \in L(\HH)}{inf}H(X)=0.$$

\vspace{0,3cm}
We now state conditions which are equivalent to the existence of minimum of \eqref{eq22}.

\begin{thm} \label{Teo3} Let $A, B \in CR(\HH), \ C \in L(\HH)$ and $W \in L(\HH)^{+}.$ 
Then the following conditions are equivalent:
\begin{enumerate}
	\item [i)] The set $\{H(X): X \in L(\HH)\}$ has a minimum, i.e., there exists $X_0 \in L(\HH)$ such that
	\begin{equation} 
	H(X) \geq H(X_0), \mbox{ for every } X \in L(\HH), \label{eq3}
	\end{equation}
	\item[ ii)]$R(C) \subseteq R(A) + R(A)^{\perp_{W}}$ and $N(B) \subseteq N(A^*WC),$
	\item [iii)] the normal equation 
	\begin{equation}
	A^*W(AXB-C)=0, \label{normaleq1}
	\end{equation}
	admits a solution.
\end{enumerate}

If any of these conditions holds, then 
$$\underset{X \in L(\HH)}{min} H(X)=C^{*}W_{/R(A)}C.$$
Moreover, the operator $X_0 \in L(\HH)$ satisfies 
$$\underset{X \in L(\HH)}{min} H(X)=H(X_0),$$
if and only if $X_0B$ is a $W$- inverse of $A$ in $R(C).$ 
\end{thm}

\begin{dem} 
\noindent \dit{$i) \Rightarrow ii)$} 
Suppose $H(X)$ has a minimum element. Let $X_0 \in L(\HH)$ such that 
$$H(X_0) \leq H(X), \mbox{ for every } X \in L(\HH),$$ or equivalently 
$$\NWW{(AX_0B-C)x} \leq \NWW{(AXB-C)x}, \mbox{ for every } X \in L(\HH) \mbox{ and } x \in \HH.$$
If $x \not \in N(B)$ then $y=Bx \not = 0,$ and given $z \in \HH,$ there exists $X \in L(\HH)$ such that $z=Xy.$ Therefore
$$\NWW{AX_0Bx-Cx} \leq \NWW{Az-Cx}, \mbox{ for every } z \in \HH.$$
Then $X_0Bx$ is a $W$-$LSS$ of $Az=Cx,$ and by Theorem \ref{thmWLSS}, $Cx \in R(A) + R(A)^{\perp_{W}},$ concluding that
$$C(\HH \setminus N(B)) \subseteq  R(A) + R(A)^{\perp_{W}}.$$

Observe that since $\HH \setminus N(B)$ is a non-empty open set, and
$$\HH \setminus N(B) \subseteq C^{-1}( R(A) + R(A)^{\perp_{W}}),$$
the subspace $C^{-1}(R(A) + R(A)^{\perp_{W}})$ has a non-empty interior, therefore 
$$\HH = C^{-1}(R(A) + R(A)^{\perp_{W}}),$$ then,
$$R(C) \subseteq R(A) + R(A)^{\perp_{W}}.$$

Observe also, that since the interior of the subspace $N(B) \subsetneq \HH $ is empty, then the set $\HH \setminus N(B)$ is a dense subset of $\HH.$ Therefore given $y \in R(C)$ there exists 
$x \in \HH$ such that $y=Cx,$ and there exists a sequence $\{x_n\}_{n \geq 1} \subset \HH \setminus N(B)$ such that $\underset{n \ra \infty}{lim \ x_n}=x.$
Then 
$$\NWW{AX_0Bx_n-Cx_n} \leq \NWW{Az-Cx_n}, \mbox{ for every } z \in \HH, \mbox{ and for every } n \in \mathbb{N},$$ and taking limit on both sides of the inequality, we get
$$\NWW{AX_0Bx-Cx} \leq \NWW{Az-Cx}, \mbox{ for every } x, z \in \HH.$$ Therefore by Theorem \ref{thmWinversa},
$G=X_0B$ is a $W$-inverse of $A$ in $R(C)$ such that $GP_{N(B)}=X_0BP_{N(B)}=0.$ Then by Theorem \ref{thmWinversa},  $A^*WC=A^*WAG$ and multiplying by $P_{N(B)}$ we get
$$A^*WCP_{N(B)} =A^{*}WAGP_{N(B)}=0,$$ and then $N(B) \subseteq N(A^*WC).$

\noindent \dit{$ii) \Rightarrow iii)$}  If $R(C) \subseteq R(A) + R(A)^{\perp_{W}}$, by Theorem \ref{thmWinversa}, there exists
$X_0 \in L(\HH)$ a solution of the normal equation
\begin{equation}
A^*W(AX_0-C)=0. \label{eqnor1}
\end{equation}  
Since $N(B) \subseteq N(A^*WC),$ we have that $A^*WC=A^*WCP_{N(B)}+ A^*WCP_{N(B)^{\perp}}=A^*WCP_{N(B)^{\perp}},$ then multiplying \eqref{eqnor1} by $P_{N(B)^{\perp}}$ it follows that
$$A^*W(AX_0P_{N(B)^{\perp}}-CP_{N(B)^{\perp}})=A^*W(A(X_0B^{\dagger})B-C)=0,$$ 
and then equation \eqref{normaleq1}, admits a solution.

\noindent \dit{$iii) \Rightarrow i)$} Let $X_0$ be a solution of the normal equation \eqref{normaleq1}, then by Theorem \ref{thmWinversa}, $G_0=X_0B$  is a  $W$-inverse of $A$ in $R(C),$ then we have 
$$\NWW{AG_0x-Cx} \leq \NWW{Az-Cx}, \mbox{ for every } x, z \in \HH.$$
Given $Y \in L(\HH),$ take $z=Yx,$ therefore
$$\NWW{AX_0Bx-Cx} \leq \NWW{AYx-Cx}, \mbox{ for every } Y \in L(\HH), \mbox{ and every } x \in \HH.$$
In particular, if $Y=XB,$ then 

$$\NWW{(AX_0B-C)x} \leq \NWW{(AXB-C)x}, \mbox{ for every } X \in L(\HH), \mbox{ and every } x \in \HH.$$
And $$H(X_0) \leq H(X), \mbox{ for every } X \in L(\HH).$$

Finally,  $X_0$ is the minimum of Problem \ref{eq22}, if and only if $X_0$ is a solution of the equation \eqref{normaleq1}, if and only if $X_0B$ is a $W$-inverse of $A$ in $R(C)$ (see Theorem \ref{thmWinversa}).
Therefore, in this case
$$H(X_0)=(AX_0B-C)^*W(AX_0B-C)=C^*W_{/R(A)}C,$$  where we used \cite[Theo.~4.3]{Contino}.
Then
$$\underset{X \in L(\HH)}{min} H(X)=H(X_0)=C^*W_{/R(A)}C.$$ 
\end{dem}

\begin{cor} \label{Cor4} Let $A, B \in CR(\HH), \ C \in L(\HH)$ and $W \in L(\HH)^{+}.$ Suppose that $R(C) \subseteq R(A) + R(A)^{\perp_{W}}$ and $N(B) \subseteq N(A^*WC),$ then the solutions of problem \eqref{eq3} (or equation \eqref{normaleq1}) are
$$(A^{*}WA)^{\dagger}A^{*}WCB^{\dagger}+L-(A^{*}WA)^{\dagger}A^{*}WALBB^{\dagger},$$ for arbitrary $L \in L(\HH).$
\end{cor}

\begin{dem} Since $R(C) \subseteq R(A) + R(A)^{\perp_{W}}$ and $N(B) \subseteq N(A^*WC),$ by Theorem \ref{Teo3}, problem \eqref{eq3} (or equation \eqref{normaleq1}) admits a solution. Then, by Theorem \ref{Teosol}, we get to the conclusion.
\end{dem}

\bigskip
\section{Minimization results in $S_p$}

In this section we study the approximation problem presented in the introduction: given $A, B \in CR(\HH)$, $C \in L(\HH)$ and $W \in L(\HH)^{+}$ such that $W^{1/2} \in S_p$ for some $p$  with $1 \leq p < \infty,$  analize the existence of
\begin{equation}
\underset{X \in L(\HH)}{min} \Nphissw{AXB-C}, \label{eq5}
\end{equation}
where $\Nphissw{X}=\Nphiss{W^{1/2}X}.$

Observe that, from equation \eqref{eqintro} and Proposition \ref{Prop Nashed}, it follows that 
$$\underset{X \in L(\HH)}{inf} \Nphissw{AXB-C} \geq \Vert W_{/ R(A)}^{1/2}C \Vert_{p}.$$

The next proposition gives sufficient conditions for the existence of minimum of \eqref{eq5}.

\begin{prop} \label{Prop3}  Let $A, B \in CR(\HH), \ C \in L(\HH)$ and $W \in L(\HH)^{+},$ such that $W^{1/2} \in S_p,$ for some $p$ with $1 \leq p < \infty.$ If
$$N(B) \subseteq N(A^*WC) \mbox{ and } R(C) \subseteq R(A) + R(A)^{\perp_{W}},$$ 
then there exists $X_0 \in L(\HH)$ such that
$$\underset{X \in L(\HH)}{min} \Vert AXB-C \Vert_{p,W}=\Vert AX_0B-C \Vert_{p,W}=\Vert W_{/R(A)}^{1/2}C \Vert_{p,W}.$$
\end{prop}

\begin{dem}
If $N(B) \subseteq N(A^*WC)$ and $R(C) \subseteq R(A) + R(A)^{\perp_{W}},$ by Theorem \ref{Teo3}, there exists $X_0\in L(\HH)$ such that $H(X_0)=\underset{X \in L(\HH)}{min}H(X)=C^{*}W_{/R(A)}C,$ i.e. $$H(X_0)=C^{*}W_{/R(A)}C\leq H(X), \mbox{ for every } X\in L(\HH).$$ 
Since $W^{1/2} \in S_p,$ by Proposition \ref{Prop Nashed}, it holds that  $$\Nphiss {W_{/R(A)}^{1/2}C}=\Nphiss{W^{1/2}(AX_0B-C)}=\Nphissw{AX_0B-C}\leq \Nphissw{AXB-C}, \mbox{ for every } X\in L(\HH),$$
then $$\underset{X \in L(\HH)}{min} \Nphissw {AXB-C} = \Nphissw{AX_0B-C}=\Nphiss {W_{/R(A)}^{1/2}C}.$$  
\end{dem}

The following result characterizes the set where the minimum of $\Nphissw{AXB-C}$ is achieve as the solutions of an equation. For the proof, we follow similar ideas as in \cite[Theo.~1.4]{Drago} and \cite[Theo.~2.6]{Mah3}.

\begin{lema}\label{critical_point}
Let $A, B \in CR(\HH),$ $C \in L(\HH)$ and $W \in L(\HH)^{+},$ such that $W^{1/2} \in S_p$ for some $p$ with $1 < p < \infty$ and  consider $F_p(X)=\Nphissw{AXB-C}^p$. Then, $X_0\in L(\HH)$ is a global minimum of $F_p$ if and only if $X_0\in L(\HH)$ is a solution of 
\begin{equation}
B\vert W^{1/2}(AXB-C) \vert^{p-1} U^{*} W^{1/2}A=0, \label{eqminimo}
\end{equation}
where $W^{1/2}(AXB-C)=U\vert W^{1/2}(AXB-C)\vert$ is the polar decomposition of the operator $W^{1/2}(AXB-C),$ with $U$ a partial isometry with $N(U)=N(W^{1/2}(AXB-C)).$
\end{lema}

\begin{dem} First observe that in \eqref{eqminimo}, $U$ varies with $X.$ 
	
Suppose $X_0$ is a global minimum of $F_p.$ Let $W^{1/2}(AX_0B-C)=U\vert W^{1/2}(AX_0B-C)\vert$ be the polar decomposition of the operator $W^{1/2}(AX_0B-C),$ with $U$ a partial isometry with $N(U)=N(W^{1/2}(AX_0B-C)).$ 
By Theorem \ref{TeoD}, $F_p$ has a $\phi-directional$ derivative for all $\phi \in [0, 2\pi).$ Then it is easy to check that, for every $\ X, \ Y\in L(\HH)$ and $\phi \in [0, 2\pi),$
$$D_{\phi}F_p (X,Y)=D_{\phi}G_p (W^{1/2}(AXB-C),W^{1/2}AYB),$$
where $G_p(X)=\Nphis{X}^{p}.$
Then, by Theorem \ref{TeoD} and Lemma \ref{LemaM2}, it holds for every $\phi \in [0, 2\pi)$
$$0\leq D_{\phi}F_p(X_0,Y)=p \ Re \ [e^{i\phi} tr ( \vert W^{1/2}(AX_0B-C) \vert^{p-1} U^{*} W^{1/2}AYB)], \mbox{ for every }Y \in L(\HH).$$
Considering a suitable $\phi$ and $Y,$ it follows that
$$B\vert W^{1/2}(AX_0B-C) \vert^{p-1} U^{*} W^{1/2}A=0.$$

Conversely, suppose that $X_0 \in L(\HH)$ is a solution of \eqref{eqminimo}, then for any $\phi \in [0,2\pi)$ and $Y \in L(\HH)$ we have
$$D_{\phi}F_p(X_0,Y)=0.$$
If $F_p(X_0)=0$ then $X_0$ is a minimum of $F_p.$  Suppose that $F_p(X_0) \not = 0$ and let $f_p(X)=F_p(X)^{\frac{1}{p}},$ then
$$D_{\phi}f_p (X_0,Y)=0, \mbox{ for every } Y \in L(\HH).$$
Let $g_p(X)=\Vert X \Vert_{p},$ then it is easy to check that, for every $Y \in L(\HH),$ we have
\begin{align*}
0&=D_{\phi}f_p(X_0,e^{i(\pi-\phi)}(-Y+X_0))\\
&=D_{\phi}g_p (W^{1/2}(AX_0B-C),W^{1/2}Ae^{i(\pi-\phi)}(-Y+X_0)B+e^{i(\pi-\phi)}W^{1/2}C-e^{i(\pi-\phi)}W^{1/2}C)\\
&=D_{\phi}g_p (W^{1/2}(AX_0B-C),-e^{i(\pi-\phi)}W^{1/2}(AYB-C)+e^{i(\pi-\phi)}W^{1/2}(AX_0B-C)).
\end{align*}

On the other hand, by \cite[Theo.~1.4]{Drago}, if $X \in L(\HH),$ we have that $$D_{\phi}g_p(X,e^{i(\pi-\phi)}X)=-\Vert X \Vert_{p}.$$ 

Then
\begin{align*}
\Vert W^{1/2}(AX_0B-C) \Vert_{p}&=-D_{\phi}g_p(W^{1/2}(AX_0B-C),e^{i(\pi-\phi)}W^{1/2}(AX_0B-C))+\\
&+D_{\phi}g_p (W^{1/2}(AX_0B-C),-e^{i(\pi-\phi)}W^{1/2}(AYB-C)+e^{i(\pi-\phi)}W^{1/2}(AX_0B-C))\\
&\leq -D_{\phi}g_p(W^{1/2}(AX_0B-C),e^{i(\pi-\phi)}W^{1/2}(AX_0B-C))+\\
&+D_{\phi}g_p (W^{1/2}(AX_0B-C), e^{i(\pi-\phi)}W^{1/2}(AX_0B-C))\\
&+D_{\phi}g_p (W^{1/2}(AX_0B-C),-e^{i(\pi-\phi)}W^{1/2}(AYB-C)) \\
&=D_{\phi}g_p (W^{1/2}(AX_0B-C),-e^{i(\pi-\phi)}W^{1/2}(AYB-C)) \\
& \leq \Vert-e^{i(\pi-\phi)} W^{1/2}(AYB-C) \Vert_{p}=\Vert W^{1/2}(AYB-C) \Vert_{p}, \mbox{ for every } Y \in L(\HH),
\end{align*}
where we used properties of $D_{\phi}g_p$ (see the definition of $D_{\phi}g_p$ ). Then $X_0$ is a  global minimum of $f_p$ or equivalently $X_0$ is a global minimum of $F_p.$
\end{dem}

\begin{thm} \label{Teo1} Let $A, B \in CR(\HH),$ $C \in L(\HH)$ and $W \in L(\HH)^{+},$ such that $W^{1/2} \in S_p$ for some $p$ with $1 \leq p < \infty$ and $N(B)\subseteq N(A^*WC).$  Then the following are equivalent:

\begin{enumerate}
\item [i)] There exists $X_0 \in L(\HH)$ such that
$$\underset{X \in L(\HH)}{min} \Nphissw{AXB-C}=\Nphissw{AX_0B-C},$$ 
\item [ii)] the normal equation 
\begin{equation}
A^{*}W(AXB-C)=0, \label{normal1}
\end{equation}
admits a solution.
\item [iii)] $R(C) \subseteq R(A)+R(A)^{\perp_{W}},$
\item [iv)] 
there exists $X_0 \in L(\HH)$ such that
$$\underset{X \in L(\HH)}{min} (AXB-C)^*W(AXB-C)=(AX_0B-C)^*W(AX_0B-C).$$
\end{enumerate}
In this case, $$\underset{X \in L(\HH)}{min} \Nphissw{AXB-C}=\Nphiss{W_{/R(A)}^{1/2}C}.$$
Moreover, $X_0 \in L(\HH)$ satisfies $$\Nphissw{AX_0B-C}=\Nphiss{W_{/R(A)}^{1/2}C},$$ if and only if $X_0$ is as in Corollary \ref{Cor4}.
\end{thm}

\begin{dem}
\noindent \dit{$i) \Rightarrow ii)$}  For $1\leq p < \infty,$ consider  
$F_p: S_p \rightarrow \mathbb{R}^{+},$ $$F_p(X)=\Nphiss{W^{1/2}(AXB-C)}^{p}.$$
By Theorem \ref{TeoD}, $F_p$ has a $\phi-directional$ derivative for all $\phi \in [0, 2\pi).$ Then it is easy to check that, for every $\ X, \ Y\in L(\HH)$ and $\phi \in [0, 2\pi),$
$$D_{\phi}F_p (X,Y)=D_{\phi}G_p (W^{1/2}(AXB-C),W^{1/2}AYB),$$
where $G_p(X)=\Nphis{X}^{p}.$

Suppose that there exists $X_0 \in L(\HH),$ a global minimum  of $\Nphissw{AXB-C}.$ Then $X_0$ is a global minimum of $F_p$ and, by Lemma \ref{LemaM2}, we have
$$\underset{0 \leq \phi < 2\pi}{inf} (D_{\phi}F_p(X_0,Y)) \geq 0, \mbox{ for every } Y \in L(\HH).$$
Let $W^{1/2}(AX_0B-C)=U\vert W^{1/2}(AX_0B-C)\vert$ be the polar decomposition of the operator $W^{1/2}(AX_0B-C),$ with $U$ a partial isometry with $N(U)=N(W^{1/2}(AX_0B-C)),$  $P=P_{N(W^{1/2}(AX_0B-C))}$ and $Q=P_{N((W^{1/2}(AX_0B-C))^{*})}.$

\vspace{0,3cm}
If $p=1,$ by Theorem \ref{TeoD} and Lemma \ref{LemaM2} it holds, for every $\phi \in [0, 2\pi)$
$$0\leq D_{\phi}F_1(X_0,Y)=Re \ [e^{i\phi} tr(U^{*}W^{1/2}AYB)] + \Vert QW^{1/2}AYBP\Vert_{1},  \mbox{ for every } Y \in L(\HH).$$
Considering a suitabe $\phi$ for each $Y \in L(\HH),$ we get
$$\vert tr(U^{*}W^{1/2}AYB) \vert \leq \Vert QW^{1/2}AYBP\Vert_{1}, \mbox{ for every } Y \in L(\HH),$$ or equivalently
$$\vert tr(U^{*}W^{1/2}AZ) \vert \leq \Vert QW^{1/2}AZP\Vert_{1}, \mbox{ for every } Z \in L(\HH), \mbox{ with } N(B) \subseteq N(Z).$$
Observe that $R(Q)=N(U^{*})$ and $R(P)=N(U),$ therefore $U^{*}Q=PU^{*}=0.$ Also, observe that since $N(B) \subseteq N(A^*WC),$ we have that 
$$N(B) \subseteq N(A^*W(AX_0B-C)(I-P)).$$ 

Let $Y \in L(\HH)$ then 
$$\vert tr(U^{*}W^{1/2}AYA^*W(AX_0B-C)) \vert  = \vert tr((I-P)U^* W^{1/2}AYA^*W(AX_0B-C)) \vert =$$
$$= \vert tr(U^* W^{1/2}AYA^*W(AX_0B-C)(I-P)) \vert \leq \Vert QW^{1/2}AYA^*W(AX_0B-C)(I-P)P\Vert_{1}=0,$$
where we used $N(B) \subseteq N(YA^*W(AX_0B-C)(I-P)).$
Then
$$tr(U^{*}W^{1/2}AYA^*W(AX_0B-C))=tr(A^*W(AX_0B-C)U^{*}W^{1/2}AY)=0, \mbox{ for every } Y \in L(\HH).$$
Therefore 
$$A^*W(AX_0B-C)U^{*}W^{1/2}A=A^*W^{1/2}U\vert W^{1/2}(AX_0B-C)\vert U^{*}W^{1/2}A= 0.$$
Hence
$$\vert W^{1/2}(AX_0B-C) \vert^{1/2} U^{*}W^{1/2}A=0,$$
or
$$\vert W^{1/2}(AX_0B-C) \vert U^{*}W^{1/2}A=(AX_0B-C)^{*}WA=0,$$
or  equivalently
$$A^*W(AX_0B-C)=0.$$
\vspace{0,3cm}

If $1 < p < \infty,$ by Lemma \ref{critical_point},
$$B\vert W^{1/2}(AX_0B-C) \vert^{p-1} U^{*} W^{1/2}A=0.$$
Observe that, since $N(B) \subseteq N(A^*WC),$ we have that
\begin{equation}
N(B) \subseteq N(A^*W(AX_0B-C))=N(A^*W^{1/2}U\vert W^{1/2}(AX_0B-C)\vert). \label{eqprueba}
\end{equation}
On the other hand we have
$$R(\vert W^{1/2}(AX_0B-C) \vert^{p-1} U^{*} W^{1/2}A) \subseteq N(B),$$ and from \eqref{eqprueba} we have that
$$A^*W^{1/2}U\vert W^{1/2}(AX_0B-C)\vert\vert W^{1/2}(AX_0B-C) \vert^{p-1} U^{*} W^{1/2}A=$$
$$A^*W^{1/2}U\vert W^{1/2}(AX_0B-C) \vert^{p} U^{*} W^{1/2}A=0.$$

Then
$$A^*W^{1/2}U\vert W^{1/2}(AX_0B-C) \vert^{p} U^{*} W^{1/2}A=$$
$$=A^*W^{1/2}U\vert W^{1/2}(AX_0B-C) \vert^{p/2}\vert W^{1/2}(AX_0B-C) \vert^{p/2} U^{*} W^{1/2}A=0.$$ Therefore

$$\vert W^{1/2}(AX_0B-C) \vert^{p/2} U^{*} W^{1/2}A=0,$$
and since $N(\vert W^{1/2}(AX_0B-C) \vert^{r})=N(\vert W^{1/2}(AX_0B-C) \vert^{s}),$ for $s, t > 0.$
We have that
$$\vert W^{1/2}(AX_0B-C) \vert U^{*} W^{1/2}A=(AX_0B-C)^*WA=0,$$ or equivalently
$$A^{*}W(AX_0B-C)=0.$$

\noindent \dit{$ii) \Rightarrow iii)$} See Theorem \ref{thmWinversa}.

\noindent \dit{$iii) \Rightarrow iv)$} It follows form Theorem \ref{Teo3}.

\noindent \dit{$iv) \Rightarrow i)$} See the proof of Proposition \ref{Prop3}.

\vspace{0,3cm}
Finally, $X_0 \in L(\HH)$ is such that
$\underset{X \in L(\HH)}{min} \Nphissw{AXB-C}=\Nphissw{AX_0B-C},$ if and only if $X_0$ is a solution of the normal equation \eqref{normal1}, and then $X_0$ is as in Corollary \ref{Cor4} and by Theorem \ref{Teo3} and Proposition \ref{Prop Nashed} $$\underset{X \in L(\HH)}{min} \Nphissw{AXB-C}=\Nphissw{AX_0B-C}=\Nphiss{W_{/R(A)}^{1/2}C}.$$
\end{dem}

\vspace{0,3cm}
Observe that the equation $A^*W(AXB-C)=0$ admits a solution, if and only if the equation $A^*W(AXB-C)B^*=0$ admits a solution and $N(B)\subseteq N(A^*WC).$ Then when equation $A^*W(AXB-C)=0$ admits a solution, the set of solutions of equation $A^*W(AXB-C)=0$ and equation $A^*W(AXB-C)B^*=0$ coincides.
Observe also, that if $N(B)\subseteq N(A^*WC),$ then $R(C) \subseteq R(A)+R(A)^{\perp_{W}}$ if and only if $R(CB^*) \subseteq R(A)+R(A)^{\perp_{W}}.$

\vspace{0,3cm}
When $p=2,$ it is possible to characterize the existence of minimum of Problem \ref{eq5}, without additional assumptions.
	
\begin{thm} \label{thm11} Let $A, B \in CR(\HH), \ C \in L(\HH)$ and $W \in L(\HH)^{+},$ such that $W^{1/2} \in S_2.$ 
Then the following are equivalent:
\begin{itemize}
	\item [i)] There exists the minimum of problem \eqref{eq5} for $p=2$, i.e., there exists $X_0 \in L(\HH)$ such that
	$$\underset{X \in L(\HH)}{min} \Vert AXB-C \Vert_{2,W}=\Vert AX_0B-C \Vert_{2,W},$$
	\item [ii)] the normal equation 
	\begin{equation}
	A^{*}W(AXB-C)B^{*}=0, \label{normaleq2}
	\end{equation}
	admits a solution.
	\item [iii)] $R(CB^*)\subseteq R(A) + R(A)^{\perp_{W}}.$
\end{itemize}
In this case, $$\underset{X \in L(\HH)}{min} \Vert AXB-C \Vert_{2,W}=\Vert W_{/R(A)}^{1/2}C \Vert_{2}.$$
Moreover, $X_0 \in L(\HH)$ satisfies $$\Vert AX_0B-C \Vert_{2,W}=\Vert W_{/R(A)}^{1/2}C \Vert_{2},$$ if and only if $X_0$ is as in Corollary \ref{Cor4}.
\end{thm}

\begin{dem}
\noi \dit{$i) \Leftrightarrow ii)$} It follows from Lemma \ref{critical_point}.

\noi \dit{$ii) \Leftrightarrow iii)$}  $R(CP_{N(B)^{\perp}})=R(CB^*)\subseteq R(A) + R(A)^{\perp_{W}}$ and $N(B) \subseteq N(CP_{N(B)^{\perp}})$ if and only if (by Theorem \ref{Teo3}), there exists a solution of the equation 
\begin{equation}
A^{*}W(AXB-CP_{N(B)^{\perp}})=0, \label{eqp2}
\end{equation}
if and only if, there exists a solution of the equation
$$A^{*}W(AXB-C)B^*=0.$$  

\vspace{0,3cm}
Finally, $X_0 \in L(\HH)$ is the minimum of Problem \ref{eq5} for $p=2,$ if and only if $X_0$ is a solution of the normal equation \eqref{normaleq2} (or equivalently $X_0$ is a solution of equation \eqref{eqp2}), then $X_0$ is as in Corollary \ref{Cor4} and by Theorem \ref{Teo3} and Proposition \ref{Prop Nashed} $$\underset{X \in L(\HH)}{min} \Vert AXB-C \Vert_{2,W}=\Vert AX_0B-C \Vert_{2,W}=\Vert W_{/R(A)}^{1/2}C \Vert_{2}.$$
\end{dem} 

The existence of solutions of \eqref{eq22} implies the existence of solutions of \eqref{eq5}, Example \ref{example3} shows that the converse it is not true, notice that $N(B) \not \subseteq N(A^*WC)$, then \eqref{eq22} has not minimum. 
%
%
%
Also this example shows that in general, for $1<p<\infty$ a global minimum of $F_p: S_p \rightarrow \mathbb{R},$ $F_p(X)=\Nphissw{AXB-C}^{p}$ is not necessarily a solution of the normal equation $A^{*}W(AXB-C)B^*=0$, which contradicts \cite[Theorem~4.1]{Mah4}.

\example \label{example3}  Let $\HH=\mathbb{C}^{2},$ $W=I$, the identity matrix, $A=\matriz{1}{0}{-1}{0},$ $B=\matriz{a^2}{-1}{a^2}{-1}$ and $C=\matriz{1}{0}{0}{a^\frac{2}{p-1}}\matriz{-1}{0}{0}{1}$, with $a,p > 1$.

Let $X_0=\matriz{1}{-1}{0}{0}$, then it is easy to verify that $AX_0B=0$, thus
$$B\vert AX_0B-C \vert^{p-1}U^*A=B\vert C \vert^{p-1} U^*A= B\matriz{-1}{0}{0}{a^2}A=0,$$ i.e., in virtue of Lemma \ref{critical_point} $X_0,$ is a global minimum of $F_p$. 

On the other hand, $$B(AX_0B-C)^*A=-BC^*A^*=B\matriz{-1}{0}{0}{a^\frac{2}{p-1}}A=\matriz{-a^2+a^\frac{2}{p-1}}{0}{-a^2+a^\frac{2}{p-1}}{0}\neq 0,$$
for every $p\neq2$. 
Then for $p\neq2,$  it follows that $X_0$ is a global minimum of $F_p$ but is not a solution of the normal equation $A^{*}W(AXB-C)B^*=0.$

\section*{Acknowledgements}
Maximiliano Contino was supported by Peruilh fundation and CONICET PIP 0168. Juan I. Giribet was partially supported by CONICET PIP 0168 and UBACyT2014.  A. Maestripieri was partially supported by CONICET PIP 0168.

\end{document}